\documentclass[12pt]{article}
\usepackage{amsmath}
\usepackage{amssymb}

\begin{document}

\begin{titlepage}
\title{\bf Lifts of Time Dependent Complex  \\ Hamiltonian Mechanical Systems}
\author{ Mehmet Tekkoyun \footnote{tekkoyun@pamukkale.edu.tr} \\
 {\small Department of Mathematics, Pamukkale University,}\\
{\small 20070 Denizli, Turkey}}
\date{\today}
\maketitle

\begin{abstract}

In this study, firstly, the k-th order extension of complex
product manifold is considered. Then the higher order vertical, complete lifts
 of geometric structures on product manifold to its extended spaces are given. Also higher
order lifts of tensor field of type (1,1) are presented. And then extended
contact manifolds are defined. Finally higher order vertical
and complete lifts of time dependent complex Hamiltonian equations
on contact manifold to its extensions are introduced. In conclusion,
geometric meaning of Hamiltonian mechanical systems is discussed.

{\bf Keywords:} complex, Kaehlerian and contact manifolds,
extensions, lifting theory, Hamiltonian mechanics.

{\bf PACS:} 02.40

\end{abstract}
\end{titlepage}

\section{\textbf{Introduction}}

Lifting theory was introduced by Bowman in 1970 \cite{bowman}. It is well
known that it permits to extend the differentiable structures. Therefore it
has an important role in differential geometry. Really, in before studies,
extensions of any real, complex manifold and complex product manifold were
defined and the higher order vertical, complete and horizontal lifts of
functions, vector fields and 1-forms on any manifold to its extension spaces
were studied in \cite{tekkoyun1, tekkoyun2, tekkoyun4, tekkoyun5} and there
in.

Modern differential geometry provides a fundamental framework for
studying Hamiltonian mechanics. In recent years, it is possible to
find many studies about differential geometric methods in
mechanics \cite{carinena, mcrampin, deleon1, deleon2} and there
in. We know that the dynamics of Hamiltonian formalisms is
characterized by a suitable vector field defined on cotangent
bundles which are phase-spaces of momentum of a given
configuration manifold. $H:T^{\ast }Q\rightarrow \mathbf{R}$ is a
regular Hamiltonian function then there is a unique vector field
$Z_{H}$ on cotangent bundle $T^{\ast }Q$ such that dynamical
equations
\begin{equation}
i_{Z_{H}}\Phi =dH,  \label{1.1}
\end{equation}%
where $\Phi $ is the symplectic form and $H$ stands for Hamiltonian
function. The paths of the Hamiltonian vector field $Z_{H}$ are the
solutions of the Hamiltonian equations shown by

\begin{equation}
\frac{dq^{i}}{dt}=\frac{\partial H}{\partial p_{i}},\,\frac{dp_{i}}{dt}=-%
\frac{\partial H}{\partial q^{i}},  \label{1.2}
\end{equation}

where $q^{i}$ and $(q^{i},p_{i}),1\leq i\leq m,$ are coordinates of $Q$ and $%
T^{*}Q.$ The triple, either $(T^{*}Q,\Phi ,Z_{H})$ or $(T^{*}Q,\Phi ,H),$ is
called \textit{Hamiltonian system} on the cotangent bundle $T^{*}Q$ with
symplectic form $\Phi $. Time dependent complex analogous of Hamiltonian
equations given in (\ref{1.2}) is the equations

\begin{equation}
\frac{dz_{i}}{dt}=\frac{1}{\mathbf{i}}\frac{\partial H}{\partial \overline{z}%
_{i}},\frac{d\overline{z}_{i}}{dt}=-\frac{1}{\mathbf{i}}\frac{\partial H}{%
\partial z_{i}}  \label{1.3}
\end{equation}

which is introduced in \cite{tekkoyun3}.

The paper is organized as follows. In section 2, we recall the k-th order
extension $^{k}N$ of a 2m+1-dimensional product manifold $N$ and the higher
order vertical, complete of functions, vector fields and 1-forms on $N$ to $%
^{k}N.$ Also, we will give the higher order vertical and complete lifts of
complex tensor field of type (1,1) on $N$ to $^{k}N$ and extended contact
manifolds structured in \cite{tekkoyun5}. In sections 3 and 4 we introduce
the higher order vertical and complete lifts of time dependent complex
Hamiltonian equations and discuss geometric results about Hamiltonian
formalisms on higher order mechanical systems.

The manifolds, tensor fields, and geometric objects we consider in this
paper, are assumed to be differentiable of class $C^{\infty }$ (i.e.,
smooth) and the sum is taken over repeated indices. Unless otherwise stated
it will be accepted $0\leq r\leq k,1\leq i\leq m.$ Also, $v$ and $c$ will
denote the vertical and complete lifts of geometric structures on either $%
^{k-1}M$ to $^{k}M$ or $^{k-1}N$ to $^{k}N.$ Dots mean derivation with
respect to time. The symbol $\mathrm{C}_{j}^{r}$ called combination is the
binomial coefficient $\binom{r}{j}.$

\section{\textbf{Preliminaries}}

In this section, we will summarize the studies given in \cite{tekkoyun5},
i.e., we recall $k$-th order extension of a complex product manifold and
higher order vertical and complete lifts of differential elements on complex
product manifold to its extension spaces. Also we present the manifolds
named to be extended contact manifolds.

\subsection{\textbf{Extended Complex Product Manifolds }}

We say to be \textit{extended complex product manifold} to the k-th order
extension $^{k}N=\mathbf{R}\times ^{k}M\ $ of 2m+1- dimensional product
manifold $N=M\times \mathbf{R},$ where $^{k}M$ extended complex manifold.
Let $(t,z^{ri},\,\overline{z}^{ri})$ be a coordinate system on a
neighborhood $^{k}V$ of any point $p$ of $^{k}N$. Therefore, by $\left\{
\newline
\frac{\partial }{\partial t},\frac{\partial }{\partial z^{ri}},\frac{%
\partial }{\partial \overline{z}^{ri}}\right\} $ and $\left\{ dt,dz^{ri},d%
\overline{z}^{ri}\right\} $ we define natural bases over coordinate system
of tangent space $T_{p}(^{k}N)$ and cotangent space $T_{p}^{\ast }(^{k}N)$
of $^{k}N,$ respectively$.$

Let $f$ be a complex function defined on $N$ and $(t,z^{0i},\overline{z}%
^{0i})$ be coordinates of $N.$ Therefore, 1- form defined by equality
\begin{equation}
df=\frac{\partial f}{\partial t}dt+\frac{\partial f}{\partial z^{0i}}dz^{0i}+%
\frac{\partial f}{\partial \overline{z}^{0i}}d\overline{z}^{0i}  \label{2.1}
\end{equation}%
is differential of $f$. Let $\chi (N)$ set of vector fields and $\chi ^{\ast
}(N)$ the set of dual vector fields on $N.$ In this case, any elements $Z$
and $\omega $ of $\chi (N)$ and $\chi ^{\ast }(N)$ are respectively
determined by
\begin{equation}
Z=\frac{\partial }{\partial t}+Z^{0i}\frac{\partial }{\partial z^{0i}}+%
\overline{Z}^{0i}\frac{\partial }{\partial \overline{z}^{0i}},  \label{2.2}
\end{equation}%
and
\begin{equation}
\omega =dt+Z_{0i}dz^{0i}+\overline{Z}_{0i}d\overline{z}^{0i},  \label{2.3}
\end{equation}%
such that $Z^{0i},\omega _{0i}\in \digamma (M).$

\subsection{\textbf{Higher Order Lifts of Geometrical Structures}}

In this section, we recall extensions of some definitions and properties
about the higher order vertical and complete of geometrical elements defined
on 2m+1- dimensional product manifold $N$ to its extension $^{k}N$. The
\textit{vertical lift} of function $f$ to $^{k}N$ is called the function $%
f^{v^{k}}$ defined by equality
\begin{equation}
f^{v^{k}}=f\circ \tau _{N}\circ \tau _{^{2}N}\circ ...\circ \tau _{^{k-1}N}.
\label{3.2}
\end{equation}%
such that a natural projection $\tau _{^{k-1}N}:^{k}N\rightarrow ^{k-1}N$ .
The \textit{complete lift} of function $f$ to $^{k}N$ is said to be the
function $f^{c^{k}}$ defined by equality
\begin{equation}
f^{c^{k}}=t(\frac{\partial f^{c^{k-1}}}{\partial t})^{v}+\overset{.}{z}^{ri}(%
\frac{\partial f^{c^{k-1}}}{\partial z^{ri}})^{v}+\overset{.}{\overline{z}}%
^{ri}(\frac{\partial f^{c^{k-1}}}{\partial \overline{z}^{ri}})^{v}.
\label{3.6}
\end{equation}%
Let $f^{c^{r}}$ be $r$ -th order complete lift of a function $f\in \digamma
(N)$ to $^{r}N.$

The \textit{vertical \ }and\textit{\ complete lifts }of vector field $Z$ on $%
N$ to $^{k}N$ are the vector field $Z^{v^{k}}$ on $^{k}N$ defined by
equality
\begin{equation}
Z^{v^{k}}(f^{c^{k}})=(Zf)^{v^{k}},\text{ }Z^{c^{k}}(f^{c^{k}})=(Zf)^{c^{k}}.
\label{4.2}
\end{equation}
Given by (\ref{2.2}) the vector field $Z$ defined on $N.$ Then, vertical \
and complete lifts of $Z$ to $^{k}N$ are
\begin{equation}
Z^{v^{k}}=\frac{\partial }{\partial t}+(Z^{0i})^{v^{k}}\frac{\partial }{%
\partial z^{ki}}+(\overline{Z}^{0i})^{v^{k}}\frac{\partial }{\partial
\overline{z}^{ki}},\text{ }Z^{c^{k}}=\frac{\partial }{\partial t}+\mathrm{C}%
_{r}^{k}(Z^{0i})^{v^{k-r}c^{r}}\frac{\partial }{\partial z^{ri}}+\mathrm{C}%
_{r}^{k}(\overline{Z}^{0i})^{v^{k-r}c^{r}}\frac{\partial }{\partial
\overline{z}^{ri}}.  \label{4.3}
\end{equation}

The \textit{vertical \ }and\textit{\ complete lifts} of 1-form $\omega $ on $%
N$ to $^{k}N$ are the 1-form $\omega ^{v^{k}}$ on $^{k}N$ defined by
equality
\begin{equation}
\omega ^{v^{k}}(Z^{c^{k}})=(\omega Z)^{v^{k}},\text{ }\omega
^{c^{k}}(Z^{c^{k}})=(\omega Z)^{c^{k}}\text{ }  \label{5.2}
\end{equation}
Denote by (\ref{2.3}) the 1-form $\omega $ defined on $N.$ Then vertical \
and complete lifts of order $k$ to $^{k}N$ of $\omega $ are
\begin{equation}
\omega ^{v^{k}}=dt+(\omega _{0i})^{v^{k}}dz^{0i}+(\overline{\omega }%
_{0i})^{v^{k}}d\overline{z}^{0i},\text{ }\omega ^{c^{k}}=dt+(\omega
_{0i})^{c^{k-r}v^{r}}dz^{ri}+(\overline{\omega }_{0i})^{c^{k-r}v^{r}}d%
\overline{z}^{ri}.  \label{5.3}
\end{equation}%
The \textit{vertical lift} of a tensor field of type (1,1) $\phi $ to $^{k}N$
is the structure $\phi ^{v^{k}}$ on $^{k}N$ given by
\begin{equation}
\phi ^{v^{k}}(\xi ^{c^{k}})=(\phi \xi )^{v^{k}},\text{ }\eta ^{v^{k}}(\phi
^{v^{k}})=(\eta \phi )^{v^{k}}.  \label{6.2}
\end{equation}%
The \textit{complete lift }of a tensor field of type (1,1) $\phi $ to $^{k}N$
is the structure $\phi ^{c^{k}}$ on $^{k}N$ given by
\begin{equation}
\phi ^{c^{k}}(\xi ^{c^{k}})=(\phi \xi )^{v^{k}},\text{ }\eta ^{c^{k}}(\phi
^{c^{k}})=(\eta \phi )^{c^{k}}.  \label{6.4}
\end{equation}

Let $^{k}N=\mathbf{R}\times ^{k}M\ $ be k-th extension of 2m+1- dimensional
product manifold $N=\mathbf{R\times }M,$ i.e., let $^{k}N$ be extended
complex product manifold$.$ A triple $(\phi ^{c^{k}},\xi ^{c^{k}},\eta
^{v^{k}})($or $(\phi ^{c^{k}},\xi ^{c^{k}},\eta ^{c^{k}}))$, is called an
\textit{extended contact structure }on $^{k}N$ such that $\phi ^{v^{k}},\phi
^{c^{k}}$ are tensors of type (1,1), $\xi ^{c^{k}}$ is a vector field and $%
\eta ^{v^{k}},\eta ^{c^{k}}$ are differential 1-forms on $^{k}N$ defined by
\begin{equation}
\phi ^{v^{k}}=-I+\xi ^{c^{k}}\otimes \eta ^{v^{k}},\,\text{\ }\eta
^{v^{k}}(\xi ^{c^{k}})=1\,{,}\text{ }\phi ^{c^{k}}=-I+\xi ^{c^{k}}\otimes
\eta ^{c^{k}},\text{ }\eta ^{c^{k}}(\xi ^{c^{k}})=1  \label{7.1}
\end{equation}%
An extended manifold $^{k}N$ endowed with a contact structure $(\phi
^{c^{k}},\xi ^{c^{k}},\eta ^{v^{k}})($or $(\phi ^{c^{k}},\xi ^{c^{k}},\eta
^{c^{k}}))$ is said to be an\textit{\ extended contact manifold}.\thinspace
\thinspace It is well-known that if $k=0,$ the manifold $N=M\times R$ with
contact structure $(\phi ,\xi ,\eta )$ is named \textit{contact manifold}.
By means of (\ref{7.1}), the higher order vertical and complete lifts of a
tensor field of type (1,1) on a contact manifold $N$ obey the following
generic properties
\begin{equation*}
\begin{array}{ll}
i) & \phi ^{v^{k}}(\xi ^{c^{k}})=0,\text{ }\phi ^{c^{k}}(\xi ^{c^{k}})=0, \\
ii) & \eta ^{v^{k}}(\phi ^{v^{k}})=0,\text{ }\eta ^{c^{k}}(\phi ^{c^{k}})=0,%
\end{array}%
\end{equation*}%
for all $\xi \in \chi (N),\eta \in \chi ^{\ast }(N)$ and $\phi \in \Im
_{1}^{1}(N),rank\phi ^{v^{k}}($or $(\phi ^{c^{k}}))=m(k+1).$\thinspace
\thinspace \thinspace \thinspace \thinspace

\section{Higher Order \textbf{Vertical Lifts of Time Dependent Hamiltonians}}

In this section, we introduce higher order vertical lifts of time dependent
Hamiltonian equations for classical mechanics structured on contact
manifold. Let $^{k}N$ be k-th order extension of contact manifold $N$ fixed
with extended coordinates $(t,z_{ri},\overline{z}_{ri})$. Then we define
vector fields $\left\{ \newline
\frac{{\partial }}{{\partial t}},\newline
\frac{{\partial }}{{\partial z^{ri}}},\newline
\frac{{\partial }}{{\partial \overline{z}^{ri}}}\right\} $ and dual covector
fields $\left\{ dt,dz_{ri},d\overline{z}_{ri}\right\} $ being bases of
tangent space $T_{p}(^{k}N)$ and cotangent space $T_{p}^{\ast }(^{k}N)$ of $%
^{k}N.$ $(\phi ^{\ast })^{c^{k}}$ and $\lambda ^{v^{k}}=(\phi ^{\ast
}{}^{c^{k}}(\omega ^{v^{k}}))$ are respectively k-th order complete and
vertical lifts of contact structure $\phi ^{\ast }$ \thinspace being the
dual of $\,\phi $ and Liouville form $\lambda $ on $N.$ If $\Phi
^{v^{k}}=-d\lambda ^{v^{k}}$ is k-th order vertical lift of closed 2-form $%
\Phi =-d\lambda ,$ then we say that $\Phi ^{v^{k}}$ is a closed 2-form on $%
^{k}M$.

\textit{A time -dependent vector field} on an extended Kaehlerian manifold $%
^{k}M$ is a $C^{\infty }$ map $Z^{v^{k}}:\mathbf{R}\times ^{k}M\rightarrow
T(^{k}M)$ such that $Z^{v^{k}}(t,p)\in T_{p}(^{k}M).$ All the results
obtained on extended Kaehlerian manifold\ $^{k}M$ hold for time dependent
vector fields. Hence, we set $\Phi _{t,s}^{v^{k}}(p)$ to be the integral
curve of $Z_{t}^{v^{k}}$ trough time $t=s$, i.e.,
\begin{equation}
\frac{d}{dt}(\phi _{t,s}^{v^{k}}(p))=Z_{t}^{v^{k}}(\phi _{t,s}^{v^{k}}(p))
\label{8.1}
\end{equation}%
and
\begin{equation}
\phi _{t,s}^{v^{k}}(p)=p,\text{ }t=s,  \label{8.2}
\end{equation}%
where $Z_{t}^{v^{k}}$ is the vector field on $^{k}M$ given by $%
Z_{t}^{v^{k}}(p)=Z^{v^{k}}(t,p).$ In fact, $\phi _{t,s}^{v^{k}}$ is the (%
\textit{time -dependent}) \textit{local 1-parameter group \thinspace
generated} by $Z_{t}^{v^{k}}.$

Let $H^{v^{k}}:^{k}N=\mathbf{R}\times ^{k}M\rightarrow \mathbf{C}$ be a
function on $^{k}N.$ For each $t\in \mathbf{R}$ we define $%
H_{t}^{v^{k}}:^{k}M\rightarrow \mathbf{C}$ by $%
H_{t}^{v^{k}}(p)=H^{v^{k}}(t,p).$ By \textit{time dependent} \textit{%
Hamiltonian vector field} we call the vector field $Z_{H_{t}}^{v^{k}}$ on $%
^{k}M$ with energy $H_{t}^{v^{k}}$ given by the isomorphism
\begin{equation}
i_{Z_{t}^{v^{k}}}\Phi ^{v^{k}}=dH_{t}^{v^{k}}.  \label{8.3}
\end{equation}%
where for simplicity, we set $Z_{t}^{v^{k}}=Z_{H_{t}^{v^{k}}}^{v^{k}}.$
Consider a mapping $Z^{v^{k}}:^{k}N=\mathbf{R}\times ^{k}M\rightarrow
T(^{k}M)$ by $Z^{v^{k}}(t,p)=Z_{t}^{v^{k}}(p)\in T_{p}(^{k}M),$ $t\in
\mathbf{R,}$ $p\in ^{k}M.$ Then there is a vector field $%
Z_{H^{v^{k}}}^{v^{k}}$ on extended contact manifold $^{k}N$ given by $%
Z_{H^{v^{k}}}^{v^{k}}(t,p)=\frac{\partial }{\partial t}+Z^{v^{k}}(t,p),$
i.e.,
\begin{equation}
Z_{H^{v^{k}}}^{v^{k}}(t,p)=\frac{\partial }{\partial t}+Z_{t}^{v^{k}}(p).
\label{8.4}
\end{equation}

\textbf{Proposition 1: }Let $^{k}M$\textbf{\ }be extended Kaehlerian
manifold with\textbf{\ }closed 2-form $\Phi ^{v^{k}}.$ The k-th order
vertical lift of Hamiltonian vector field $Z_{t}$ on Kaehlerian manifold $M$
endowed with\textbf{\ }closed 2-form $\Phi $ is given by
\begin{equation}
Z_{t}^{v^{k}}=\frac{1}{{\mathbf{i}}}\frac{\partial H^{v^{k}}}{\partial
\overline{z}_{0i}}\frac{\partial }{\partial z_{ki}}-\frac{1}{{\mathbf{i}}}%
\frac{\partial H^{v^{k}}}{\partial z_{0i}}\frac{\partial }{\partial
\overline{z}_{ki}}.  \label{8.5}
\end{equation}

\textbf{Proof:} Let $^{k}M$\textbf{\ }be extended Kaehlerian manifold with%
\textbf{\ }closed 2-form $\Phi ^{v^{k}}.$ Consider that $Z_{t}^{v^{k}}$ is
the k-th order vertical lift of Hamiltonian vector field $Z_{t}$ associated
Hamiltonian energy $H_{t}.$ Also, $Z_{t}^{v^{k}}$ is Hamiltonian vector
field $Z_{t}^{v^{k}}$ associated Hamiltonian energy $H_{t}^{v^{k}}$and given
by
\begin{equation}
Z_{t}^{v^{k}}=(Z_{0i})^{v^{k}}\frac{\partial }{\partial z_{ri}}+(\overline{Z}%
_{0i})^{v^{k}}\frac{\partial }{\partial \overline{z}_{ri}}  \label{8.6}
\end{equation}
For the closed 2-form $\Phi ^{v^{k}}$ on $^{k}M,$ we find
\begin{equation}
\Phi ^{v^{k}}=-d\lambda ^{v^{k}}=-d(\frac{1}{2}{\mathbf{i}}(-z_{0i}d%
\overline{z}_{0i}+\overline{z}_{0i}dz_{0i}))=-{\mathbf{i}}d\overline{z}%
_{0i}\wedge dz_{0i}.  \label{8.7}
\end{equation}%
Taking into consideration the isomorphism given in (\ref{8.3}), we calculate
by
\begin{equation}
i_{Z_{t}^{v^{k}}}\Phi ^{v^{k}}=i_{Z_{t}^{v^{k}}}(-d\lambda ^{v^{k}})=\mathbf{%
i}(\overline{Z}_{0i})^{v^{k}}dz_{0i}+\mathbf{i}(Z_{0i})^{v^{k}}d\overline{z}%
_{0i}.  \label{8.8}
\end{equation}%
On the other hand, the differential of Hamiltonian energy $H_{t}^{v^{k}}$ on
$^{k}M$ we give by
\begin{equation}
dH_{t}^{v^{k}}=\frac{\partial H_{t}^{v^{k}}}{\partial z_{ri}}dz_{ri}+\frac{%
\partial H_{t}^{v^{k}}}{\partial \overline{z}_{ri}}d\overline{z}_{ri}.
\label{8.9}
\end{equation}%
From equality (\ref{8.3})$,$ the k-th order vertical lift of Hamiltonian
vector field on Kaehlerian manifold $M$ fixed with\textbf{\ }closed 2-form $%
\Phi $ we find as
\begin{eqnarray*}
Z_{t}^{v^{k}} &=&\frac{1}{{\mathbf{i}}}\frac{\partial H_{t}^{v^{k}}}{%
\partial \overline{z}_{ri}}\frac{\partial }{\partial z_{ri}}-\frac{1}{%
\mathbf{i}}\frac{\partial H_{t}^{v^{k}}}{\partial z_{ri}}\frac{\partial }{%
\partial \overline{z}_{ri}} \\
&=&\frac{1}{{\mathbf{i}}}\frac{\partial H_{{}}^{v^{k}}}{\partial \overline{z}%
_{ri}}\frac{\partial }{\partial z_{ri}}-\frac{1}{{\mathbf{i}}}\frac{\partial
H_{{}}^{v^{k}}}{\partial z_{ri}}\frac{\partial }{\partial \overline{z}_{ri}}.
\end{eqnarray*}
Thus, proof is complete. $\Box $

Suppose that the curve
\begin{equation}
\alpha ^{k+1}:I=(-\epsilon ,\epsilon )\subset \mathbf{C}\rightarrow \mathbf{R%
}\times ^{k}M=^{k}N  \label{8.10}
\end{equation}%
be an integral curve of Hamiltonian vector field $Z_{H}^{v^{k}},$with $%
\epsilon >0,$ i.e.,
\begin{equation}
\overset{.}{\alpha }^{k+1}(t)=Z_{H}^{v^{k}}(\alpha (t)),\,t\in I.
\label{8.11}
\end{equation}%
In the local coordinates we have
\begin{equation}
\alpha ^{k+1}(t)=(t,z_{ri}(t),\overline{z}_{ri}(t)).  \label{8.12}
\end{equation}%
So, we obtain
\begin{equation}
\overset{.}{\alpha }^{k+1}(t)=Z_{H}^{v^{k}}(t,z_{ri}(t),\overline{z}%
_{ri}(t))=\frac{\partial }{\partial t}+Z_{t}^{v^{k}}(z_{ri}(t),\overline{z}%
_{ri}(t)).  \label{8.13}
\end{equation}%
Now, from $\overset{.}{\alpha }^{k+1}(t)=Z_{H}^{v^{k}}(\alpha ^{k+1}(t)),$
then we infer the following equations
\begin{equation}
\frac{dz_{ri}}{dt}=\frac{1}{{\mathbf{i}}}\frac{\partial H^{v^{k}}}{\partial
\overline{z}_{0i}},\frac{d\overline{z}_{ri}}{dt}=-\frac{1}{{\mathbf{i}}}%
\frac{\partial H^{v^{k}}}{\partial z_{0i}},  \label{8.14}
\end{equation}%
that is called \textit{k-th order} \textit{vertical lift of time dependent
complex Hamiltonian equations} on contact manifold $N$.

In (\ref{8.14}), if k=0, we get the equations
\begin{equation}
\frac{dz_{0i}}{dt}=\frac{1}{{\mathbf{i}}}\frac{\partial H}{\partial
\overline{z}_{0i}},\frac{d\overline{z}_{0i}}{dt}=-\frac{1}{{\mathbf{i}}}%
\frac{\partial H}{\partial z_{0i}},  \label{8.15}
\end{equation}%
or
\begin{equation}
\frac{dz_{i}}{dt}=\frac{1}{{\mathbf{i}}}\frac{\partial H}{\partial \overline{%
z}_{i}},\frac{d\overline{z}_{i}}{dt}=-\frac{1}{{\mathbf{i}}}\frac{\partial H%
}{\partial z_{i}},  \label{8.16}
\end{equation}%
which is \textit{time dependent complex Hamiltonian equations} on contact
manifold $N$ determined in (\ref{1.3}) and introduced in \cite{tekkoyun3}.

\section{Higher Order \textbf{Complete Lifts of Time Dependent Hamiltonians}}

In this section, we bring in higher order complete lifts of time dependent
complex Hamiltonian equations for classical mechanics structured on contact
manifold. Let $^{k}N$ be k-th order extension of contact manifold $N$ and
endowed with extended coordinates $(t,z_{ri},\overline{z}_{ri})$. Then by $%
\left\{ \frac{\partial }{\partial t},\frac{\partial }{\partial z_{ri}},\frac{%
\partial }{\partial \overline{z}_{ri}}\right\} $ and $\left\{ dt,dz_{ri},d%
\overline{z}_{ri}\right\} ,$ we determine vector fields and dual covector
fields being bases of tangent space $T_{p}(^{k}N)$ and cotangent space $%
T_{p}^{\ast }(^{k}N)$ of $^{k}N,$ respectively$.$ We define by $(\phi ^{\ast
})^{c^{k}}$ and $\lambda ^{c^{k}}=(\phi ^{\ast }{}^{c^{k}}(\omega ^{c^{k}}))$
the k-th order complete lifts of contact structure $\phi ^{\ast }$
\thinspace being the dual of $\,\phi $ and Liouville form $\lambda $ on $M,$
respectively$.$ If $\Phi ^{c^{k}}=-d\lambda ^{c^{k}}$ is k-th order vertical
lift of closed 2-form $\Phi =-d\lambda ,$ then we say that $\Phi ^{c^{k}}$
is a closed 2-form on extended Kaehlerian manifold $^{k}M$. \textit{A time
-dependent vector field} on an extended Kaehlerian manifold $^{k}M$ is a $%
C^{\infty }$ map $Z^{c^{k}}:^{k}N\rightarrow T(^{k}M)$ such that $%
Z^{c^{k}}(t,p)\in T_{p}(^{k}M).$ All the results obtained on extended
Kaehlerian manifold\ $^{k}M$ hold for time dependent vector fields. Hence,
we set $\Phi _{t,s}^{c^{k}}(p)$ to be the integral curve of $Z_{t}^{c^{k}}$
trough time $t=s$, i.e.,
\begin{equation}
\frac{d}{dt}(\phi _{t,s}^{c^{k}}(p))=Z_{t}^{c^{k}}(\phi _{t,s}^{c^{k}}(p))
\label{9.1}
\end{equation}%
and
\begin{equation}
\phi _{t,s}^{c^{k}}(p)=p,\text{ }t=s,  \label{9.2}
\end{equation}%
where $Z_{t}^{c^{k}}$ is the vector field on $^{k}M$ given by $%
Z_{t}^{c^{k}}(p)=Z^{c^{k}}(t,p).$ In fact, $\phi _{t,s}^{c^{k}}$ is the
\textit{(time -dependent ) local 1-parameter group \thinspace generated} by $%
Z_{t}^{c^{k}}.$

Let $H^{c^{k}}:^{k}N\rightarrow \mathbf{C}$ be a function on $^{k}N.$ For
each $t\in \mathbf{R}$ we define $H_{t}^{c^{k}}:^{k}M\rightarrow \mathbf{C}$
by $H_{t}^{c^{k}}(p)=H^{c^{k}}(t,p).$ By \textit{time dependent} \textit{%
Hamiltonian vector field} we say to be the vector field $Z_{H_{t}}^{c^{k}}$
on $^{k}M$ with energy $H_{t}^{c^{k}}$ given by the isomorphism
\begin{equation}
i_{Z_{t}^{c^{k}}}\Phi ^{c^{k}}=dH_{t}^{c^{k}}.  \label{9.3}
\end{equation}%
where for simplicity, we set $Z_{t}^{c^{k}}=Z_{H_{t}^{c^{k}}}^{c^{k}}.$
Define a mapping $Z^{c^{k}}:^{k}N\rightarrow T(^{k}M)$ by $%
Z^{c^{k}}(t,p)=Z_{t}^{c^{k}}(p)\in T_{p}(^{k}M),$ $t\in \mathbf{R,}$ $p\in
^{k}M.$ Then there is a vector field $Z_{H^{c^{k}}}^{c^{k}}$ on extended
contact manifold $^{k}N$ given by $Z_{H^{c^{k}}}^{c^{k}}(t,p)=\frac{\partial
}{\partial t}+Z^{c^{k}}(t,p),$ i.e.,
\begin{equation}
Z_{H^{c^{k}}}^{c^{k}}(t,p)=\frac{\partial }{\partial t}+Z_{t}^{c^{k}}(p).
\label{9.4}
\end{equation}%
\textbf{Proposition 2: }Let $^{k}M$\textbf{\ }be extended Kaehlerian
manifold with\textbf{\ }closed 2-form $\Phi ^{c^{k}}.$ The k-th order
complete lift of Hamiltonian vector field $Z_{t}$ on Kaehlerian manifold $M$
fixed with\textbf{\ }closed 2-form $\Phi $ is given by
\begin{equation}
Z_{t}^{c^{k}}=\frac{1}{{\mathbf{i}}}\frac{\partial H^{c^{k}}}{\partial
\overline{z}_{ri}}\frac{\partial }{\partial z_{ri}}-\frac{1}{{\mathbf{i}}}%
\frac{\partial H^{c^{k}}}{\partial z_{ri}}\frac{\partial }{\partial
\overline{z}_{ri}}.  \label{9.5}
\end{equation}

\textbf{Proof:} Let $^{k}M$\textbf{\ }be extended Kaehlerian manifold with%
\textbf{\ }closed Kaehlerian form $\Phi ^{c^{k}}$. Consider that $%
Z_{t}^{c^{k}}$ be the k-th order complete lift of Hamiltonian vector field $%
Z_{t}$ associated Hamiltonian energy $H_{t}$ and given by
\begin{equation}
Z_{t}^{c^{k}}=\mathrm{C}_{r}^{k}(Z_{0i})^{v^{k-r}c^{r}}\frac{\partial }{%
\partial z_{ri}}+\mathrm{C}_{r}^{k}(\overline{Z}_{0i})^{v^{k-r}c^{r}}\frac{%
\partial }{\partial \overline{z}_{ri}}.  \label{9.6}
\end{equation}%
For the closed Kaehlerian form $\Phi ^{c^{k}}$ on $^{k}M,$ we obtain
\begin{equation}
\Phi ^{c^{k}}=-d\lambda ^{c^{k}}=-d(\frac{1}{2}{\mathbf{i}}(-z_{ri}d%
\overline{z}_{ri}+\overline{z}_{ri}dz_{ri}))=-{\mathbf{i}}d\overline{z}%
_{ri}\wedge dz_{ri}.  \label{9.7}
\end{equation}%
Using by the isomorphism given in (\ref{9.3}), we find
\begin{equation}
\begin{array}{l}
i_{Z_{t}^{c^{k}}}\Phi ^{c^{k}}=-{\mathbf{i}}\mathrm{C}_{r}^{k}(\overline{Z}%
_{0i})^{v^{k-r}c^{r}}dz_{ri}+{\mathbf{i}}\mathrm{C}%
_{r}^{k}(Z_{0i})^{v^{k-r}c^{r}}d\overline{z}_{ri}%
\end{array}
\label{9.8}
\end{equation}%
On the other hand, the differential of Hamiltonian energy $H_{t}^{c^{k}}$we
define by
\begin{equation}
dH_{t}^{c^{k}}=\frac{\partial H_{t}^{c^{k}}}{\partial z_{ri}}dz_{ri}+\frac{%
\partial H_{t}^{c^{k}}}{\partial \overline{z}_{ri}}d\overline{z}_{ri}.
\label{9.9}
\end{equation}%
By means of equality (\ref{9.3})$,$ the Hamiltonian vector field $%
Z_{t}^{c^{k}}$on extended Kaehlerian manifold $^{k}M$ is calculated as
follows:
\begin{eqnarray*}
Z_{t}^{c^{k}} &=&\frac{1}{{\mathbf{i}}}\frac{\partial H_{t}^{c^{k}}}{%
\partial \overline{z}_{ri}}\frac{\partial }{\partial z_{ri}}-\frac{1}{%
\mathbf{i}}\frac{\partial H_{t}^{c^{k}}}{\partial z_{ri}}\frac{\partial }{%
\partial \overline{z}_{ri}} \\
&=&\frac{1}{{\mathbf{i}}}\frac{\partial H_{{}}^{c^{k}}}{\partial \overline{z}%
_{ri}}\frac{\partial }{\partial z_{ri}}-\frac{1}{{\mathbf{i}}}\frac{\partial
H_{{}}^{c^{k}}}{\partial z_{ri}}\frac{\partial }{\partial \overline{z}_{ri}}.
\end{eqnarray*}
Hence, proof finishes. $\Box $

Assume that the curve
\begin{equation}
\alpha ^{k+1}:I=(-\epsilon ,\epsilon )\subset \mathbf{C}\rightarrow \mathbf{R%
}\times ^{k}M=^{k}N  \label{9.10}
\end{equation}%
be an integral curve of Hamiltonian vector field $Z_{H}^{c^{k}},$with $%
\epsilon >0,$ i.e.,
\begin{equation}
\overset{.}{\alpha }^{k+1}(t)=Z_{H}^{c^{k}}(\alpha (t)),\,t\in I.
\label{9.11}
\end{equation}%
In the local coordinates it holds
\begin{equation}
\alpha ^{k+1}(t)=(t,z_{ri}(t),\overline{z}_{ri}(t)).  \label{9.12}
\end{equation}%
Therefore we have
\begin{equation}
\overset{.}{\alpha }^{k+1}(t)=Z_{H}^{c^{k}}(t,z_{ri}(t),\overline{z}%
_{ri}(t))=\frac{\partial }{\partial t}+Z_{t}^{c^{k}}(z_{ri}(t),\overline{z}%
_{ri}(t)).  \label{9.13}
\end{equation}%
Now, from $\overset{.}{\alpha }^{k+1}(t)=Z_{H}^{c^{k}}(\alpha ^{k+1}(t)),$
then the equations obtained by%
\begin{equation}
\frac{dz_{ri}}{dt}=\frac{1}{{\mathbf{i}}}\frac{\partial H^{c^{k}}}{\partial
\overline{z}_{ri}},\frac{d\overline{z}_{ri}}{dt}=-\frac{1}{\mathbf{i}}\frac{%
\partial H^{c^{k}}}{\partial z_{ri}},  \label{9.14}
\end{equation}
are \textit{k-th order complete lift of time dependent complex Hamiltonian
equations} on contact manifold $N$.

In (\ref{9.14}), if k=0, we have the equations
\begin{equation}
\frac{dz_{0i}}{dt}=\frac{1}{{\mathbf{i}}}\frac{\partial H}{\partial
\overline{z}_{0i}},\frac{d\overline{z}_{0i}}{dt}=-\frac{1}{{\mathbf{i}}}%
\frac{\partial H}{\partial z_{0i}},  \label{9.15}
\end{equation}%
or
\begin{equation}
\frac{dz_{i}}{dt}=\frac{1}{{\mathbf{i}}}\frac{\partial H}{\partial \overline{%
z}_{i}},\frac{d\overline{z}_{i}}{dt}=-\frac{1}{{\mathbf{i}}}\frac{\partial H%
}{\partial z_{i}},  \label{9.16}
\end{equation}%
which is \textit{time dependent complex Hamiltonian equations} on contact
manifold $N$ given in (\ref{1.3}) and obtained in \cite{tekkoyun3}.

\textbf{Corollary: }By means of the equations found the above, we conclude
that the\textbf{\ }Hamiltonian formalisms in generalized classical mechanics
and field theory can be intrinsically characterized on the extended contact
manifolds $^{k}N$, and the geometric approach of complex Hamiltonian systems
is that the solutions of time dependent vector fields $Z_{t}^{v^{k}}$ and $%
Z_{t}^{c^{k}}$ on extended Kaehlerian manifolds $^{k}M$ are paths time
dependent complex Hamiltonian equations obtained (\ref{8.14}) and (\ref{9.14}%
) on extended contact manifolds $^{k}N,$ respectively$.$ Hence, by means of
the lifting theory, it is shown that Hamiltonian formalism may be
generalized to extended contact manifolds $^{k}N$.

\end{document}